\documentclass[reqno]{amsart}
\usepackage{amssymb,amsmath,amsthm,hyperref,comment,xcolor}

\newtheorem{theorem}{Theorem}

\newtheorem{lemma}[theorem]{Lemma}

\renewcommand{\leq}{\leqslant}
\renewcommand{\geq}{\geqslant}

\newcommand*{\bbe}{
  \mathop{
    \mathchoice{\vcenter{\hbox{\larger[4]$\mathbb{E}$}}}
               {\kern0pt\mathbb{E}}
               {\kern0pt\mathbb{E}}
               {\kern0pt\mathbb{E}}
  }\displaylimits
}

\newcommand{\abs}[1]{\left\lvert #1\right\rvert}

\begin{document}
\title{Optimal bounds for the length and number of unit fraction representations} 

\author[T. F. Bloom]{Thomas F. Bloom}
\address{University of Manchester, United Kingdom}
\email{bloommaths@gmail.com} 

\author[W. van Doorn]{Wouter van Doorn}
\address{Groningen, the Netherlands}
\email{wonterman1@hotmail.com}

\date{}

\begin{abstract}This paper gives a human-written exposition and simplification of proofs that were recently posted by OpenAI~\cite{AI26}. We prove that every $a/b\in (0, 1)$ can be written as the sum of $O(\log\log b)$ many distinct unit fractions. We also prove that the number of representations of $1$ as the sum of $k$ distinct unit fractions grows doubly exponentially in $k$.
\end{abstract}

\maketitle

\section{Introduction}

The study of unit fractions\footnote{Often called Egyptian fractions, after the ancient Egyptian practice of writing all rational numbers as the sum of unit fractions.} is one of the oldest areas in number theory. It was, in particular, a favourite topic of Erd\H{o}s and Graham, who asked many fascinating questions concerning unit fractions. There has been great progress made on understanding many of these questions; see, for example, \cite{BlEl22} for a recent survey covering aspects of this progress that we do not address here.

In October 2026 OpenAI released hundreds of AI-generated proofs, addressing problems from across pure mathematics. One of these papers  \cite{AI26}  concerned unit fractions, and resolved three questions which were asked by Erd\H{o}s and Graham (in \cite{ErGr80}, among other places). The proofs were formalised in Lean, and hence one can verify their correctness mechanically; yet the AI-generated accompanying paper is difficult to read, contains a lot of unnecessary material, and makes it difficult to see what the new ideas are. We therefore thought it would be useful to present our own exposition of these proofs, incorporating some simplifications. For more discussion on how the proofs in this paper differ from the AI-generated proofs, see Section~\ref{sec-AI}. 

In the remainder of this introduction we will briefly state the three main results; further historical details are given in the relevant sections. These results address problems \#304, \#148, and \#293 respectively in the online collection of Erd\H{o}s problems \cite{Bl26}.

The first (and, in our view, most interesting) result concerns how efficiently a given rational number can be represented as the sum of distinct unit fractions. 

\begin{theorem}\label{th-main}
Any rational $a/b\in (0,1)$ can be written as the sum of $O(\log\log b)$ distinct unit fractions.
\end{theorem}

The second result concerns how many ways there are to write $1$ as the sum of $k$ distinct unit fractions.

\begin{theorem}\label{th-main2}
If $F(k)$ counts how many ways $1$ can be written as the sum of $k$ distinct unit fractions then\footnote{In this paper we make frequent use of Vinogradov notation $\gg$, which means the inequality holds up to an absolute multiplicative constant. In other words, $f\ll g$ is equivalent to $f=O(g)$. We also use the convention that we interpret $\log x = 1$ whenever $x < e$, to ensure that e.g. $\log \log b$ remains strictly positive for all integers $b\geq1$.}
\[\log\log F(k) \gg k.\]
\end{theorem}

A similar bound holds for the number of ways of writing any fixed $a/b\in (0,1]$ as a sum of $k$ distinct unit fractions. The third result concerns which denominators can appear in a unit fraction representation of $1$.

\begin{theorem}\label{th-main3}
If $v(k)$ is the least integer $n\geq 2$ which does not appear as a denominator in any representation of $1$ as the sum of $k$ distinct unit fractions then
\[\log\log v(k)\gg k.\]
\end{theorem}

All three bounds are, up to the implicit constants, the best possible. The complementary bounds which show this are classical and have been known at least since \cite{Er50} if not before; we reproduce the elementary proofs in Appendix~\ref{app}. 

We also note that an easy argument shows that $v(k)\leq kF(k)+2$, and so Theorem~\ref{th-main2} is an immediate consequence of Theorem~\ref{th-main3}. We also include an alternative direct proof of Theorem~\ref{th-main2}.

In the following sections we will give proofs of these three theorems, together with further historical remarks and context. We conclude with a brief discussion of the implicit constants; we have chosen to present the proofs of the main results in a mostly constant-free fashion, to better highlight the main ideas. In Section~\ref{sec-AI} we give a brief commentary on the AI-generated paper \cite{AI26}, and compare the proofs as presented here with the AI version. To keep this paper self-contained we have also reproduced some existing arguments from the literature in Appendix~\ref{app}.

\subsection*{Acknowledgements} TB is a Royal Society University Research Fellow. ChatGPT Astra was used in the early stages to help explain the proofs in \cite{AI26}. This paper was entirely written by the authors.

\section{Efficient unit fraction representations} \label{efficient}

It has long been known that every positive rational number can be written as the sum of distinct unit fractions; this can be proved with a simple greedy algorithm, as noted by Fibonacci in the 13th century. Although this provides a finite number of summands, it is often far from the most efficient representation. For example, the greedy algorithm applied to $641/796$ writes it as the sum of $13$ distinct unit fractions, of which the largest denominator has $1348$ digits, when in fact there is a much more efficient representation in
\[\frac{641}{796}=\frac{1}{2}+\frac{1}{4}+\frac{1}{20}+\frac{1}{199}+\frac{1}{3980}\]
which uses only $5$ unit fractions.

For $a/b\in (0,1)$ the greedy algorithm terminates in at most $a$ steps, although it is conjectured (first by Schinzel in \cite{Si56}) that $3$ unit fractions already suffice as long as $b$ is sufficiently large in terms of $a$. The most notable special case of this is the Erd\H{o}s-Straus conjecture that $4/b$ can be written as the sum of $3$ distinct unit fractions for all $b \ge 3$. Here we are instead concerned with the dependence on the denominator $b$.

This question first appeared in the literature in a paper of Erd\H{o}s \cite{Er50}, and it was later popularised further by Erd\H{o}s and Graham \cite{ErGr80}. For any $1\leq a< b$ define $N(a,b)$ to be the smallest $k$ such that $a/b$ can be written as the sum of $k$ distinct unit fractions, and let $N(b)$ be the maximum of $N(a,b)$ over all $1 \leq a < b$.

The first non-trivial bound of $N(b) \ll \frac{\log b}{\log\log\log b}$ was proved by de Bruijn (as reported in \cite{Er50}, although this was never published). Erd\H{o}s \cite{Er50} provided both a lower bound as well as an improved upper bound, and showed
\[\log\log b\ll N(b) \ll \frac{\log b}{\log\log b}.\]
The lower bound is attained by $\frac{b-1}{b}$, which we will prove in Appendix~\ref{app}. The upper bound was later refined by Vose \cite{Vo85} to $N(b) \ll \sqrt{\log b},$ and AI-generated proofs \cite{vD26} recently improved this to $N(b)\ll (\log\log b)^2$. Before turning to the proof of the optimal $N(b) \ll \log\log b$, we will describe the existing strategy for obtaining upper bounds, as pioneered by Erd\H{o}s~\cite{Er50} and used in all other works on the topic since.

To begin, we note that thanks to a theorem of Takenouchi \cite{Ta21} (proved in Appendix~\ref{app}), any representation as a sum of $k$ unit fractions can be converted into a sum of $k$ distinct unit fractions. This elementary fact, often overlooked in papers on this topic, is simple to prove, but makes investigating $N(b)$ significantly simpler, as we can ignore the distinct denominators requirement.

\subsection{Practical numbers}

A positive integer $n$ is $k$-practical if every positive integer smaller than $n$ can be written as the sum of at most $k$ divisors of $n$. This definition is due to Srinivasan \cite{Sr48}, although Fibonacci used them in his study of unit fractions in the 13th century.

To see why these are useful in the context of unit fractions, let $\frac{a}{b} \in (0, 1)$ be given. We then choose a practical number $n > b$ and let $m$ be such that $\frac{m}{n} \leq \frac{a}{b} < \frac{m+1}{n}$. As $n$ is assumed to be $k$-practical, we can write both $m<n$ and $an - bm < b<n$ as sums of divisors of $n$, say $m = u_1 + \cdots + u_r$ and $an - bm = v_1 + \cdots + v_s$, where $r,s\leq k$. With this, we write
\begin{align*}
\frac{a}{b} &= \frac{m}{n} + \frac{an - bm}{bn} \\
&= \frac{u_1 + \cdots + u_r}{n} + \frac{v_1 + \cdots + v_s}{bn} \\
&= \frac{1}{n/u_1} + \cdots + \frac{1}{n/u_r} + \frac{1}{bn/v_1} + \cdots + \frac{1}{bn/v_s}
\end{align*}
as a sum of at most $2k$ unit fractions. (Note the $u_i$ and $v_j$ being divisors of $n$ is essential to ensure that all denominators are integers.)

This means that one strategy for obtaining good upper bounds for $N(b)$ is to identify an infinite sequence of practical numbers, in which $k$ is small as a function of $n$. For example, Erd\H{o}s \cite{Er50} used that $n!$ is $n$-practical. The improved upper bounds of Vose \cite{Vo85} and the recent AI-generated proofs \cite{vD26} follow this strategy, constructing ad hoc practical sequences which require fewer summands.

The key observation to improve on this strategy is the simple one that we did not need to write $m$, say, as the sum $m=u_1+\cdots+u_r$ where the $u_i$ are divisors of $n$; it would suffice to write
\begin{equation} \label{semipractical}
m=\frac{u_1}{w_1}+\cdots+\frac{u_r}{w_r},
\end{equation}
where the $u_i$ are still divisors of $n$, but the denominators $w_i$ can be arbitrary. In other words, it suffices to write every $m<n$ as the sum of a small number of fractions whose numerators are divisors of $n$. (One could call such $n$ semipractical.) 

Perhaps surprisingly, with this relaxed notion of semipracticality, the number of summands required for sufficiently small $m<n$ can be bounded by an absolute constant, even for simple sequences such as the factorials or the sequence of lowest common multiples. The proof writes $m$ as a sum of primes $p$ using Vinogradov's theorem, and then uses ideas from sum-product theory to find smooth integers that add up to a multiple of $p$. Taken together, they lead to an expression in the form of~\eqref{semipractical}. (Although this constitutes a fairly simple combination of known results, this idea of using finite field sum-product results to elevate additive representations using primes to additive representations using fractions with a particular structure is, we believe, both new and delightful.)

There is one slight wrinkle, which is that the number of summands required to write $m$ is only $O(1)$ if the integer $m$ is small; this means that this strategy only produces an efficient representation for $a/b$ if $a\leq \sqrt{b}$, say. For general $a<b$ we begin by subtracting unit fractions following a greedy algorithm to reduce to this case.

\subsection{Proof of Theorem~\ref{th-main}}

We start by reminding the reader of some classical concepts from analytic number theory. We write $P(y)=\prod_{p\leq y}p$ for the product of all primes at most $y$. An application of Mertens' first theorem then implies $P(y)\leq e^{O(y)}$, which is sufficient for our purposes. An integer $n$ is $y$-smooth if every prime dividing $n$ is $\leq y$, and it is squarefree if it is not divisible by any square $>1$.

We require the following sum-product result, which was proved by Hart and Iosevich \cite{HaIo08}.
\begin{lemma}\label{lem-sp}
If $A\subseteq \mathbb{F}_p$ satisfies $\abs{A}\geq p^{3/4}$ then $\mathbb{F}_p^*\subseteq AA+AA$. In particular $-1\in AA+AA$.
\end{lemma}
For our purposes weaker results would also suffice: for example, that there exists some $c>0$ such that if $\abs{A}\geq p^{1-c}$ then $-1$ can be written as the sum and product of $O(1)$ many elements from $A$. There are a number of results of this type; perhaps the first appeared in a paper of Glibichuk \cite{Gl07}. (Such results can also be quickly deduced from earlier work, such as the sum-product result of Bourgain, Katz, and Tao \cite{BKT04}, or from exponential sum estimates of Vinogradov \cite{Vi47}.)

We wish to apply Lemma~\ref{lem-sp} with $A$ being the set of squarefree $y$-smooth numbers (so that every element of $A$ divides $P(y)$), but it is hard to prove $\abs{A\pmod{p}} \geq p^{3/4}$ for all primes $p$ in the range we are considering. Fortunately, it suffices to show that this inequality holds for almost all primes, which can be established via a short elementary argument.
\begin{lemma}\label{lem-prime}
For all except $x^{3/4+o(1)}$ many primes $p\leq x$ there exist squarefree $O(\log x\log\log x)$-smooth integers $u_1,v_1,u_2,v_2$ with
\[u_1v_1+u_2v_2\equiv -1\pmod{p}.\]
\end{lemma}
\begin{proof}
Let $x$ be given and define $y=C\log x\log\log x$ for a sufficiently large constant $C$. With $A$ the set of squarefree $y$-smooth integers, we then find
\[\abs{A} = 2^{\pi(y)} \geq 2x^{3/4},\]
where the final inequality holds if $C$ is sufficiently large. By Lemma~\ref{lem-sp} it suffices to bound the number of primes $p\leq x$ such that $\abs{A\pmod{p}}<p^{3/4}$. 

Let $Q$ be the set of such primes and consider the number of triples $(a,b,p)$ where $a\neq b\in A$, $p\in Q$ and $a\equiv b\pmod {p}$. On one hand, since $0<\abs{a-b}\leq P(y)$, the integer $\abs{a-b}$ has at most $y$ distinct prime factors, so the number of such triples is trivially $\leq \abs{A}^2y$.

On the other hand, if $p\in Q$ then $A$ occupies $<p^{3/4}\leq x^{3/4}$ residue classes so there are at least 
\[\frac{\abs{A}^2}{x^{3/4}} - \abs{A} \gg \abs{A}^2/x^{3/4}\]
distinct $a,b\in A$ whose difference is divisible by $p$. Hence, the number of such triples is $\gg \abs{Q}\abs{A}^2/x^{3/4}$. Comparing the upper and lower bounds implies
\[\abs{Q}\ll x^{3/4}y = x^{3/4+o(1)}\]
as claimed.
\end{proof}

We can now prove the main lemma establishing the kind of semipracticality condition mentioned above: all small integers can be written as the sum of $O(1)$ many fractions with smooth numerators.

This problem does not appear to have been considered before; it should be compared to the problem of writing integers as the sum of few smooth numbers (in other words, insisting the denominators be $1$). For example, Balog and S\'{a}rk\"{o}zy \cite{BaSa84} proved that all $m$ can be written as the sum of $3$ integers, each of which is $\exp(O(\sqrt{\log m\log\log m}))$-smooth. At the much stronger level of smoothness we require, namely $O(\log m\log\log m)$,  no $O(1)$ bound on the number of summands is possible, since there are only $m^{o(1)}$ integers $\leq m$ which are this smooth. The flexibility of allowing arbitrary denominators is therefore essential.

\begin{lemma}\label{lem-rep}
Every integer $m\geq 1$ can be written as the sum of at most $10$ fractions of the shape $uv/w$ where $u,v,w\geq 1$ are positive integers, and $u$ and $v$ are squarefree and $O(\log m\log\log m)$-smooth.
\end{lemma}
\begin{proof}
By increasing the constants in the $O(\cdot)$ notation as necessary, and using an extra summand of $+1$ for even integers, it suffices to prove that every sufficiently large odd integer $m$ can be written as the sum of $9$ such fractions.

Call a prime $p\leq m$ \emph{good} if it satisfies the conclusion of Lemma~\ref{lem-prime}, and \emph{bad} otherwise. By Vinogradov's theorem \cite{Vi47} every large odd integer $m$ can be written as the sum of three primes in at least $\gg \frac{m^2}{(\log m)^3}$ many ways. Since there are at most $m^{3/4+o(1)}$ many bad primes, the number of such representations of $m$ involving a bad prime is $\ll m^{7/4+o(1)}$. Hence, provided $m$ is large enough, there exists a representation $m=p_1+p_2+p_3$ where all $p_i$ are good.

By the definition of a good prime there exists, for $1\leq i\leq 3$, squarefree $O(\log m\log\log m)$-smooth integers $u_j^{(i)},v_j^{(i)}$ such that
\[u_1^{(i)}v_1^{(i)} + u_2^{(i)}v_2^{(i)} + 1 = p_iw_i\]
for some $w_i\geq 1$. We then finish the proof by writing
\begin{equation*}
m = p_1 + p_2 + p_3 = \sum_{i=1}^3 \left(\frac{u_1^{(i)}v_1^{(i)}}{w_i} + \frac{u_2^{(i)}v_2^{(i)}}{w_i} + \frac{1}{w_i} \right). \qedhere
\end{equation*}
\end{proof}

Note that, at the cost of replacing the constant $10$ by some larger unspecified $O(1)$, the full strength of Vinogradov's theorem is not needed; it would suffice to know that, given any fixed set of primes with positive relative density, every large integer is the sum of $O(1)$ many primes from this set. This in turn can be established using elementary methods, such as Schnirelmann's elementary argument coupled with Brun's sieve~\cite{Sch33}.

We may now prove Theorem~\ref{th-main}. We first prove this in the special case when the numerator is small, in which case the proof follows along similar lines to the earlier strategy with practical numbers. After this we use a simple greedy algorithm to reduce to the case of small numerator, completing the proof of Theorem~\ref{th-main}.

\begin{lemma}\label{lem-aux}
If $1\leq a\leq \sqrt{b}$ then $a/b$ can be written as the sum of $O(\log\log b)$ many unit fractions.
\end{lemma}
The bound $a\leq \sqrt{b}$ is not particularly special here; the same argument would work for $a\leq b^{1-c}$ for any $c>0$, with the implicit constants depending on $c$.
\begin{proof}
Let $a/b$ with $1 \leq a \leq \sqrt{b}$ and $b$ sufficiently large be given and write $y = c \log b$. Here, $c > 0$ is a sufficiently small constant so that, by defining $N = P(y)^2$, we have $N < \sqrt{b}$. We also let $B$ be the largest integer $m$ for which the smoothness parameter $O(\log m\log\log m)$ from Lemma~\ref{lem-rep} is at most $y$. With these definitions, Lemma~\ref{lem-rep} then implies that every $m\leq B$ can be written as the sum of at most $10$ fractions, each of whose numerator divides $N$.

By definition of $B$ we get $\log B \gg \frac{\log b}{\log\log b}$. That is, there exists a $k=O(\log\log b)$ such that $b\leq B^k$. We now write
\[aNB^k = bd+e\]
where $0\leq e<b\leq B^k$. Since $aN <b$ we also have $0\leq d<B^k$. Writing $d$ and $e$ in base $B$ we therefore have, for some $d_i,e_i\in [0,B)$,
\[d=\sum_{i<k}d_iB^i\textrm{ and }e=\sum_{i<k}e_iB^i,\]
whence
\[\frac{a}{b}=\frac{d}{NB^k}+\frac{e}{bNB^k}=\sum_{i<k}\frac{d_i}{NB^{k-i}}+\sum_{i<k}\frac{e_i}{bNB^{k-i}}.\]
Since $d_i < B$, it follows from Lemma~\ref{lem-rep} that $d_i$ can be written as the sum of at most $10$ fractions $u_jv_j/w_j$ with $u_jv_j$ dividing $N$, so that 
\[\frac{d_i}{NB^{k-i}} = \sum_j \frac{u_jv_j}{w_jNB^{k-i}} = \sum_j \frac{1}{w_jNB^{k-i}/(u_jv_j)}, \]
where the $w_jNB^{k-i}/(u_jv_j)$ are all positive integers. As this holds for all $0 \leq i < k$, as well as for all the $e_i$, we conclude that $a/b$ is the sum of at most $20k\ll \log\log b$ many unit fractions.
\end{proof}

\begin{proof}[Proof of Theorem~\ref{th-main}]
By Takenouchi's theorem (Theorem~\ref{th-tak}) it suffices to write $a/b$ as the sum of $O(\log\log b)$ many unit fractions that are not necessarily distinct. Lemma~\ref{lem-aux} does the job if $a\leq \sqrt{b}$. To reduce to this case, we subtract unit fractions greedily: let $a/b=a_0/b_0$ and, for $i\geq 1$, choose $x_i\geq 2$ such that
\[\frac{1}{x_i}\leq \frac{a_{i-1}}{b_{i-1}}<\frac{1}{x_i-1}, \]
and let
\[a_i=a_{i-1}x_{i}-b_{i-1}\quad\textrm{and}\quad b_{i}=b_{i-1}x_{i},\]
so that $a_{i}/b_{i}=a_{i-1}/b_{i-1}-1/x_i$. 

Note that the choice of $x_i$ ensures that $a_i<a_{i-1}$ (so this process terminates in at most $a$ steps). We continue until either $a_r/b_r=0$, or $b_r\geq b^2$. At this point we have $a_r\leq a\leq b\leq \sqrt{b_r}$, and hence Lemma~\ref{lem-aux} represents $a_r/b_r$ with $O(\log\log b_r)$ many unit fractions.

The total number of unit fractions required for $a/b$ is 
\[\ll r+\log\log b_r.\]
To bound $b_r$ note that $x_i=\lceil b_{i-1}/a_{i-1}\rceil \leq b_{i-1}$, and so $b_i\leq b_{i-1}^2$. In particular, since $r$ is minimal such that $b_r\geq b^2$, we have $b_r\leq b_{r-1}^2 < b^4$. 

To bound $r$, note that if $y_i=a_i/b_i$ then
\[y_i=y_{i-1}-\frac{1}{x_i}<y_{i-1}-\frac{y_{i-1}}{y_{i-1}+1}<y_{i-1}^2,\]
and hence by induction $y_i < 2^{-2^{i-1}}$ for all $i\geq 1$. (For $i=1$ note that $a_1/b_1 < 1/2$ since either $a/b < 1/2$ or if $a/b \geq 1/2$ then $x_1=2$ and $a_1/b_1=a/b-1/2 < 1/2$.) If we have not halted with $y_i=0$ then $y_i\geq 1/b_i$, and hence 
\begin{equation} \label{bibound}
2^{2^{i-1}}\leq b_i
\end{equation} 
for all $i\geq 0$. Recalling the bound $b_r < b^4$ we deduce that $r \ll \log\log b_r\ll \log\log b$, and the proof is complete.
\end{proof}

\section{Counting unit fraction representations of unity}\label{sec-main2}
Let $F(k)$ be the number of solutions to
\[1=\frac{1}{x_1}+\cdots+\frac{1}{x_k}\]
where $1\leq x_1<\cdots<x_k$ are distinct integers. It is a folklore observation (recorded in the introduction to \cite{Ko14}, for example) that $F(k)$ satisfies $F(k)\leq F(k+1)$ for all $k\geq 3$; this can be seen by noting that the map on ordered tuples
\[(x_1<\cdots <x_k)\mapsto (x_1<\cdots<x_{k-1}<x_k+1<x_k(x_k+1))\]
is injective and preserves the property of reciprocals summing to $1$, since
\[\frac{1}{x}=\frac{1}{x+1}+\frac{1}{x(x+1)}.\]

It is a classical fact that $F(k)$ is finite, and in fact grows at most doubly exponentially in $k$, i.e.
\begin{equation} \label{uppercount}
\log\log F(k) \ll k.
\end{equation}
(We prove this in Appendix~\ref{app}.) There has been some work on improving the implicit constants in this upper bound; the best current bounds are proved by Elsholtz and Planitzer \cite{ElPl21}. Until recently the best lower bound available was
\[\log\log F(k)\gg \frac{k}{\log k},\]
proved by Konyagin \cite{Ko14}. A comparable lower bound for the counting problem in which we restrict to only odd denominators was proved by Elsholtz \cite{El16}.

In this section we prove Theorem~\ref{th-main2}, which gives the improved lower bound
\[\log\log F(k) \gg k.\]
By~\eqref{uppercount} this is optimal up to the implicit constants. One way to prove this lower bound is as an immediate consequence of the bound in Theorem~\ref{th-main3}, thanks to the simple observation that
\[F(k) \geq \frac{v(k) - 2}{k}.\]
Indeed, this inequality holds since, by definition of $v(k)$, all integers in $[2,v(k))$ occur in at least one of the $k$-tuples counted by $F(k)$. In the rest of this section we will give an alternative proof which couples Theorem~\ref{th-main} with the strategy of Konyagin.

\subsection{Konyagin's strategy}
The idea of Konyagin's proof is as follows. Using elementary identities he writes $1$ as the sum of $k-1$ distinct unit fractions, where one denominator $n$ has many divisors. This is useful since, if $d\mid n$ with $d<n$, then
\begin{equation}\label{eq-kon}
\frac{1}{n}=\frac{1}{n+d}+\frac{1}{n+\frac{n^2}{d}}
\end{equation}
where the denominators on the right-hand side are distinct. Inserting this representation yields a representation of $1$ as the sum of $k$ unit fractions; there may be duplication of one of the new denominators with one of the unchanged $k-2$ denominators, but this can only happen for at most $2(k-2)$ many $d$. This implies 
\[F(k)\geq \tau(n)-1-2(k-2).\]
To give a good lower bound for $F(k)$ it therefore suffices to find a single representation of $1$ as a sum of $k-1$ distinct unit fractions, in which one denominator has many divisors.

Konyagin used a series of elementary identities to find, for many $k$ and any $m\ll 3^k$, a representation of $1$ as the sum of $k-1$ unit fractions, where one denominator is divisible by $2^m+1$. He then used the elementary inequality (for odd $m$) $\tau(2^m+1) \geq 2^{\tau(m)-1}$. The proof concludes by choosing $m$ as the product of the first $O(k/\log k)$ odd primes, so that $m\ll 3^k$ as required, and $\log \tau(m) \gg k/\log k$. Combining these facts yields
\[\log\log F(k)\gg \log \log \tau(2^m+1)\gg \log \tau(m) \gg \frac{k}{\log k}.\]

A variant of this strategy (with $2^n+1$ replaced by $2^n-1$, and a different set of elementary identities) was used by Elsholtz \cite{El16} to prove a lower bound of comparable strength even if we restrict to odd denominators. 

\subsection{The new bound}
To obtain the new bound we use exactly the strategy above, except that the initial unit fraction decomposition is obtained by applying Theorem~\ref{th-main} to $\frac{Q-1}{Q}$, where $Q$ is a primorial (and hence has many divisors) -- adding $1/Q$ to this representation yields a representation of $1$ as the sum of $O(\log\log Q)$ unit fractions where one denominator has many divisors.

\begin{proof}[Proof of Theorem~\ref{th-main2}]
Let $Q$ be the product of the first $r$ odd primes, so that $\log\log Q\ll \log r$. By Theorem~\ref{th-main} we can write $\frac{Q-1}{Q}$ as the sum of $O(\log r)$ many unit fractions. By adding $1/Q$ we obtain a representation of $1$ as the sum of $k=O(\log r)$ unit fractions, in which one denominator is $Q$.

These denominators may not be distinct (as $Q$ could appear twice) but applying Theorem~\ref{th-tak} produces a representation of $1$ as the sum of $k$ distinct unit fractions in which one denominator is divisible by $Q$. By the argument surrounding \eqref{eq-kon} we deduce that 
\[F(k+1)\geq \tau(Q)-1-2(k-1)\gg 2^r,\]
since $Q$ has $2^r$ divisors. It follows that
\[\log\log F(k+1)\geq \log r-O(1)\gg k.\]
This bound does not apply for all $k$, but it does provide, for any $r\geq 1$, some $k=O(\log r)$ for which it holds. The monotonicity property $F(k)\leq F(k+1)$ completes the proof for all $k$.
\end{proof}

\section{The smallest missing denominator}

In this section we discuss Theorem~\ref{th-main3}; this states that, if $v(k)$ is the smallest integer $\geq 2$ which does not appear in any representation of $1$ as the sum of $k$ distinct unit fractions, then 
\begin{equation}\label{eq-vk}
\log\log v(k) \gg k.
\end{equation}
The study of $v(k)$ was initiated by Erd\H{o}s and Graham \cite{ErGr80}, who noted that the bound $v(k) \gg k!$ should follow from the work of Bleicher and Erd\H{o}s \cite{BlEr76}, but speculated that the stronger lower bound of \eqref{eq-vk} may hold.

The upper bounds $\log\log F(k)\ll k$ and $v(k)\leq kF(k)+2$, as noted in the previous section, imply that $\log\log v(k)\ll k$, so \eqref{eq-vk} is, up to constants, the best possible bound. The first serious investigation of $v(k)$ after \cite{ErGr80} was undertaken by van Doorn and Tang \cite{vDTa25}, who proved $\log v(k) \gg k^2$.

\begin{proof}[Proof of Theorem~\ref{th-main3}]
By \cite[Lemma 2.1]{vDTa25} it suffices to prove that for any $b$, there is a representation of $1$ with $O(\log\log b)$ many distinct unit fractions, including $b$ as a denominator. 

This is equivalent to representing $\frac{b-1}{b}$ as the sum of $O(\log\log b)$ many distinct unit fractions, none of which are $1/b$. We therefore run a similar greedy procedure as in the proof of Theorem~\ref{th-main}, where we assume $b\geq 7$ without loss of generality. Let $y_0=(b-1)/b-1/2-1/3$, so that $y_0<1/6$, and while $y_{i-1}>0$ choose $x_i\geq 7$ minimal such that $x_i\neq b$ and $1/x_i \leq y_{i-1}$, and let $y_i=y_{i-1}-1/x_i$.

Note that, by minimality,
\[\frac{1}{x_i}\leq y_{i-1}<\frac{1}{x_i-2}\]
(the $-2$ instead of $-1$ is because we may have to skip the denominator $b$). Therefore
\[y_i=y_{i-1}-\frac{1}{x_i}<y_{i-1}-\frac{y_{i-1}}{2y_{i-1}+1}<2y_{i-1}^2.\]
By induction it follows that $\log\log(1/y_i) \gg i$ for all $i\geq 1$. Similarly, since
\[\frac{1}{x_{i+1}}\leq y_i=y_{i-1}-\frac{1}{x_i} < \frac{2}{x_i(x_i-2)}\leq \frac{1}{x_i}\]
we also have $x_i< x_{i+1}$, so the denominators $x_i$ are all distinct. 

We continue this procedure for $r$ steps until either $y_i=0$ for some $i\leq r$ or $y_r<1/b$. Since $\log\log(1/y_i)\gg i$ the latter occurs for some $r=O(\log\log b)$. If some $y_i=0$ for $i\leq r$ then we are done, having represented $\frac{b-1}{b}$ as the sum of $O(\log\log b)$ distinct unit fractions, none of which are $1/b$.

Otherwise
\[\frac{b-1}{b}=\frac{1}{2}+\frac{1}{3}+\frac{1}{x_1}+\cdots+\frac{1}{x_r}+y_r\]
where $2<3<x_1<\cdots <x_r$ are distinct denominators, none of which are $b$. If $b_i= b_{i-1}x_i$ then $y_i$ is a fraction with denominator $b_i$, and
\[x_i-2< y_{i-1}^{-1}\leq b_{i-1}\]
and hence (assuming $b\geq 2$) we have $x_i \leq 2b_{i-1}$, whence $b_i \leq 2b_{i-1}^2$. It follows that
\[b_r\leq O(b)^{2^r},\]
and hence by Theorem~\ref{th-main} the remainder $y_r$ can be written as the sum of 
\[\ll \log\log b_r \ll r+\log\log b\ll \log\log b\]
distinct unit fractions. Since $y_r<1/b$ and $y_r<1/x_r<1/x_i$ for all $i\leq r$ none of the denominators in this representation of $y_r$ can be $b$ or any $x_i$, and hence we have written $\frac{b-1}{b}$ as the sum of $O(\log\log b)$ distinct unit fractions, none of which are $1/b$ as required.
\end{proof}

\section{On explicit constants}
Even though Theorem~\ref{th-main},~\ref{th-main2} and~\ref{th-main3} are all optimal up to the implied constants, that still leaves room to determine what the correct constants should be.\footnote{There is disagreement among the authors to what extent this question is interesting.}

To this end, recall that we defined $C$ in the proof of Lemma~\ref{lem-prime} such that, with $y = C \log x \log \log x$, we have $2^{\pi(y)} \ge 2x^{3/4}$. By the prime number theorem it follows that one can take $C = 3/(4 \log 2) + o(1)$. That is, the implicit constant in the smoothness parameter in Lemma~\ref{lem-prime} and~\ref{lem-rep} is $3/(4 \log 2) + o(1)$. 

As for Lemma~\ref{lem-aux}, one can apply the prime number theorem once again to take the constant $c$ in the definition of $y = c \log b$ to be $\frac{1}{4} + o(1)$. This ensures that we still have $N = P(y)^2 < \sqrt{b}$. Combining these constants we find 
\[\log B = \frac{(\log 2 + o(1)) \log b}{3 \log \log b}, \]
implying
\[k = \left(\frac{3}{\log 2} + o(1)\right) \log \log b.\]
Since the proof of Lemma~\ref{lem-aux} uses at most $20k$ unit fractions, whereas the initial greedy reduction uses another $r$ unit fractions (where $r \leq \frac{\log \log b}{\log 2} + O(1)$ by~\eqref{bibound}), the total number of required fractions is at most 
\[20k + r = \left(\frac{61}{\log 2} + o(1)\right)\log \log b.\]
The number of fractions needed for the greedy reduction cannot, in general, be improved beyond $\frac{\log \log b}{\log 2} + O(1)$ (see Appendix~\ref{app} for more context), but there are various ways to reduce the remaining at most $\left(\frac{60}{\log 2} + o(1)\right)\log \log b$ fractions that we use, so let us very briefly comment on two of these ideas. 

The first thing to note is that we only have to consider those $a$ in Lemma~\ref{lem-aux} with $a = b^{o(1)}$. One can check that this does not significantly increase the number of fractions needed in the greedy reduction, but it does allow us to increase $c$ from $\frac{1}{4} + o(1)$ to $\frac{1}{2} + o(1)$ for a factor of two improvement.

For another factor of two improvement, recall that we wrote $aNB^k$ as $bd + e$ in the proof of Lemma~\ref{lem-aux}, after which $d$ and $e$ were written in their base $B$ expansion with at most $k$ terms each. A different strategy, with analogous definitions of $B$, $N$ and $k$, starts with the greedy reduction to obtain a fraction $a_r/b_r$, after which $a_rN$ is written as the sum of at most $k$ terms $de$ with $d \mid b_r$ and $e \leq B$. Such an expansion is possible thanks to the fact that, by construction, $b_r$ contains a divisor in every interval $[x, bx]$ with $x \leq b_r$. Combining both of these ideas gives a final result that $\left(\frac{16}{\log 2} + o(1)\right)\log \log b$ unit fractions suffices.

Stronger versions of Lemma~\ref{lem-sp} or stronger results on prime numbers from multiplicative or additive number theory can bring the constant down further still, but this would significantly lessen the elementary nature of the arguments, and it is unclear if the optimal $\frac{1}{\log 2}$ is even achievable this way. One way to, in principle, reach this would be to lower the smoothness parameter from Lemma~\ref{lem-prime} to $o(\log x \log \log x)$. However, the lemma would then no longer be true as the number of such smooth integers $u_i, v_i$ is simply insufficient to cover almost all primes. Hence, different ideas might be necessary in order to obtain the optimal constant.


Finally, we note that every proof of Theorem~\ref{th-main} along similar lines that gives, up to $o(1)$ terms, an implied constant of $c$, then implies a corresponding constant of $\frac{1}{c}$ for Theorem~\ref{th-main2} and~\ref{th-main3}.

\section{Comparison to the AI-generated paper}\label{sec-AI}
The AI-generated paper \cite{AI26} proves Theorem~\ref{th-main} in a different, much more elaborate way; this proof occupies Sections 2-5 of that paper. After this, in Section 6 it proves Theorem~\ref{th-main2} using the divisor based argument we present in Section~\ref{sec-main2}. Finally, in Sections 7 and 8 it proves Theorem~\ref{th-main3}. As we have done here, this is done by proving a slightly stronger form of Theorem~\ref{th-main} in which we find a representation of $\frac{b-1}{b}$ as a sum of $O(\log\log b)$ unit fractions omitting $1/b$. 

This latter proof is done by a completely different, much more elementary, method than it used to prove Theorem~\ref{th-main}. This offers an alternative, and much simpler, route to proving Theorem~\ref{th-main}. Some elementary preparations, in particular Takenouchi's theorem and the greedy reduction process, appear in Section 2, and are still used; however the technical arguments of Sections 3-5 (occupying 15 pages, out of 33 total) are entirely redundant.

Strangely, nowhere in the paper is this oddity mentioned: the AI appears to have found the complicated proof of Theorem~\ref{th-main} first, then moved on later to prove Theorem~\ref{th-main3} with a different method, and not noticed that this actually provides a simpler proof for Theorem~\ref{th-main}. The paper also does not make the observation that, since $v(k)\leq kF(k)+2$, the lower bound for $F(k)$ immediately follows from the lower bound for $v(k)$. 

Our proof of Theorem~\ref{th-main} is extracted from the simpler proof of Theorem~\ref{th-main3} in \cite{AI26}, with some further simplifications made. (Most notably \cite{AI26} employs an extra iterative procedure which we avoid by passing to the base-$B$ representations in Lemma~\ref{lem-aux}.)

The original complicated proof of Theorem~\ref{th-main} in \cite{AI26} proceeds via a similar strategy: after some initial reductions with a greedy algorithm, it finds an efficient representation of a wide range of rationals as a sum of unit fractions. Instead of using the sum-product/sum of primes trick employed in the simpler proof (and in this paper), however, it finds such a representation using Fourier methods (closely related to work of Croot (see, for example, \cite{Cr03}) on other questions concerning unit fractions). This requires difficult estimates on the moments of restricted divisor functions and a delicate probabilistic argument. The authors of this paper have not yet read those sections in detail, as they seemed unnecessary for the results discussed here (although they may well contain ideas useful for deeper investigation into unit fractions).

Reassuringly, most of the important references to earlier work do appear in \cite{AI26}, although it is often unclear which parts of its lengthy arguments appear in these earlier references, and which are novel. This seems to be a common flaw in the current generation of AI-generated papers -- after listing some key references they launch into a lengthy proof, with no attempt made to explain which parts of the proof come from the cited references, and which parts are new. To be fair, the same criticism can be applied to many human-written papers, but it makes reading these AI-generated papers difficult, when the main question an expert has is `what are the new ideas'?

The main flaw of \cite{AI26} (aside from the fact that almost half of its length is taken up with a redundant alternative proof) is its verbosity. As frequent readers of AI-generated mathematics are accustomed to by now, almost every proof is written in more detail than is required, often introducing unnecessary notation and variables. Easy parts of the arguments, which are standard and usually omitted in research papers, are gone over in unnecessary levels of detail, and paid the same attention as the new steps of the proof on which the entire argument rests.

This makes it difficult to extract what the important ideas of the proof are, and can be a frustrating experience for an expert. It also does not help a newcomer to the area, who would be better served by first encountering these `standard' arguments in another setting before seeing how they are applied here.

This is not meant to be a wider criticism of AI-generated proofs or expositions, but rather an observation of the shortcomings of \cite{AI26}. The quality of AI writing can improve significantly after several iterative rounds of passing the AI-generated paper back to an AI, especially when given the specific criticisms mentioned above. Doubtless many or all of the simplifications and observations made in this paper could also be made by an AI when fed the paper \cite{AI26}.

Nonetheless, we found the process of writing this paper and thinking how to simplify the proofs to help in our own understanding. While AI may be able to produce \emph{some} exposition of these ideas, it could not have produced \emph{our} exposition of them, which contains many personal choices: which context is relevant, which steps are important and novel and should be expanded, and which steps are routine and should be condensed. 

The relative quality of expositions is up to the reader to decide. We hope that the choices we have made, and this exposition we have prepared, are helpful to some readers who would like to explore the fascinating world of unit fractions. 

\appendix

\section{Other results}\label{app}
In this appendix we gather a couple of elementary arguments that are already well-known in the area.

\subsection{Takenouchi's theorem}

The following was proved by Takenouchi \cite[Sections~II--III]{Ta21}. 
\begin{theorem}\label{th-tak}
If $a/b\in (0,1]$ is the sum of $k\geq 3$ unit fractions then it is the sum of $k$ distinct unit fractions.

Furthermore, if it is the sum of $k$ unit fractions, one of whose denominators is divisible by an odd number $n$, then it is the sum of $k$ distinct unit fractions, one of whose denominators is divisible by $n$.
\end{theorem}
\begin{proof}
Suppose
\[\frac{a}{b}=\frac{1}{x_1}+\cdots+\frac{1}{x_k}.\]
Note that, since $k\geq 3$ and $a/b\leq 1$, any repeated denominator must be $\geq 3$. If $x_i=x_j=2y$ for some $y\geq 2$ then we produce a new representation using the identity
\[\frac{1}{2y}+\frac{1}{2y}=\frac{1}{y+1}+\frac{1}{y(y+1)}.\]
If $x_i=x_j=2y+1$ for some $y\geq 1$ then we produce a new representation using the identity
\[\frac{1}{2y+1}+\frac{1}{2y+1}=\frac{1}{y+1}+\frac{1}{(y+1)(2y+1)}.\]
These substitutions preserve the total number of summands; they also preserve the property of at least one denominator being divisible by any fixed odd number $n$.

It remains to check that this procedure must terminate eventually; this follows from noting that the sum of the denominators strictly increases after each substitution. The sum of denominators is bounded simply because there are only finitely many ways of writing $a/b$ as the sum of $k$ unit fractions, which itself can be seen by induction on $k$, and noting that there must exist a denominator which is $\leq kb/a$.
\end{proof}

\subsection{Optimality}

In this section we give short elementary proofs that provide complementary bounds for the quantities considered in this paper. The proofs here are probably best attributed to `folklore'; certainly similar arguments appear in both \cite{Er50} and \cite{Ta21}. The main observation behind the bounds is the following.
\begin{lemma} \label{folk1}
If 
\[1=\frac{1}{x_1}+\cdots+\frac{1}{x_k}\]
with $1\leq x_1\leq\cdots \leq x_k$ then $x_k \leq k^{2^{k-1}}$.
\end{lemma}
\begin{proof}
Note that for $0\leq i<k$ the quantity $1-\sum_{j\leq i}\frac{1}{x_j}$ is strictly positive and can be written as a rational number with denominator $x_1\cdots x_i$, where we take the empty product to be equal to $1$. Hence,
\[\frac{1}{x_1\cdots x_i}\leq 1-\sum_{j\leq i}\frac{1}{x_j}=\sum_{i<j\leq k}\frac{1}{x_j}\leq \frac{k}{x_{i+1}}.\]
It follows that, for all $0\leq i<k$,
\[x_{i+1}\leq kx_1\cdots x_i.\]
By induction we then deduce that
\[x_i\leq k^{2^{i-1}}\]
for all $i$ with $1\leq i\leq k$. Taking $i=k$ completes the proof.
\end{proof}

A variant of this argument, which we leave as an illuminating exercise to the reader,\footnote{See \cite{BlEl22} for a simple proof of a slightly weaker bound. This stronger bound perhaps first appeared in the literature in a paper of Curtiss \cite{Cu22}. A more recent proof was given by Soundararajan \cite{So05}.} proves the stronger bound that 
\[x_k\leq c^{2^k}\]
where $c\approx 1.26408$ is the Vardi constant, defined as 
\[c=\lim u_n^{2^{-n}}\]
where $u_1=1$ and $u_{i+1}=u_i(u_i+1)$. 

We will now discuss a couple of applications of Lemma~\ref{folk1}. The first argument is from \cite{Er50}: recall $N(a,b)$ as defined in Section~\ref{efficient} and consider two representations
\[
\frac{a}{b} = \frac{1}{x_1} + \cdots + \frac{1}{x_k}\quad\textrm{and}\quad\frac{b-1-a}{b} = \frac{1}{y_1} + \cdots + \frac{1}{y_l},\]
with $k = N(a, b)$ and $l = N(b-a-1,b)$. By combining these and adding $\frac{1}{b}$  we obtain
\[1 = \frac{1}{x_1} + \cdots + \frac{1}{x_k} + \frac{1}{y_1} + \cdots + \frac{1}{y_l} + \frac{1}{b}.\]
By applying Lemma~\ref{folk1} we deduce $b\leq (k+l+1)^{2^{k+l}}$ and so
\[k+l=N(a, b) + N(b-1-a, b) \gg \log \log b.\]
In particular, for half of all $1\leq a<b$ we have $N(a,b)\gg \log\log b$, and so certainly we must have $N(b)\gg \log\log b$. One can be a little more explicit about for which $a$ this lower bound is attained: since $N(a, b) \leq a$ by the greedy algorithm, it follows that $N(a, b) \gg \log \log b$ for all $a$ with $b-a = o(\log \log b)$ (in particular for $a = b-1$, say). 

Secondly, Lemma~\ref{folk1} immediately implies upper bounds for $F(k)$, which counts the number of $1\leq x_1<\cdots <x_k$ such that $1=\frac{1}{x_1}+\cdots+\frac{1}{x_k}$, since Lemma~\ref{folk1} yields
\[F(k) \leq \binom{k^{2^{k-1}}}{k},\]
whence
\[\log\log F(k) \ll k.\]
The bound $v(k) \leq kF(k)+2$ yields $\log\log v(k)\ll k$ also.

\bibliographystyle{plainurl}
\bibliography{refs}

@PREAMBLE{{\providecommand{\noopsort}[1]{}}}

@article{Bl26,
 AUTHOR = {T. F. Bloom},
 TITLE = {www.erdosproblems.com},
 URL = {https://www.erdosproblems.com},
 YEAR = {accessed 2026-10-08},
}

@article {Cr03,
    AUTHOR = {Croot, III, Ernest S.},
     TITLE = {On a coloring conjecture about unit fractions},
   JOURNAL = {Ann. of Math. (2)},
  FJOURNAL = {Annals of Mathematics. Second Series},
    VOLUME = {157},
      YEAR = {2003},
    NUMBER = {2},
     PAGES = {545--556},
      ISSN = {0003-486X,1939-8980},
   MRCLASS = {11D68 (11N64)},
  MRNUMBER = {1973054},
MRREVIEWER = {K.\ Soundararajan},
       DOI = {10.4007/annals.2003.157.545},
       URL = {https://doi-org.manchester.idm.oclc.org/10.4007/annals.2003.157.545},
}

@article {Sr48,
    AUTHOR = {Srinivasan, A. K.},
     TITLE = {Practical numbers},
   JOURNAL = {Current Sci.},
  FJOURNAL = {Current Science},
    VOLUME = {17},
      YEAR = {1948},
     PAGES = {179--180},
      ISSN = {0011-3891},
   MRCLASS = {10.0X},
  MRNUMBER = {27799},
MRREVIEWER = {D.\ H.\ Lehmer},
}

@article {Si56,
    AUTHOR = {Sierpi\'nski, W.},
     TITLE = {Sur les d\'ecompositions de nombres rationnels en fractions
              primaires},
   JOURNAL = {Mathesis},
  FJOURNAL = {Mathesis. Filosof\'ia e Historia de las Ciencias
              Matem\'{a}ticas},
    VOLUME = {65},
      YEAR = {1956},
     PAGES = {16--32},
      ISSN = {0185-6200},
   MRCLASS = {10.0X},
  MRNUMBER = {78385},
MRREVIEWER = {P.\ Erd\H{o}s},
}

@article {BKT04,
    AUTHOR = {Bourgain, J. and Katz, N. and Tao, T.},
     TITLE = {A sum-product estimate in finite fields, and applications},
   JOURNAL = {Geom. Funct. Anal.},
  FJOURNAL = {Geometric and Functional Analysis},
    VOLUME = {14},
      YEAR = {2004},
    NUMBER = {1},
     PAGES = {27--57},
      ISSN = {1016-443X,1420-8970},
   MRCLASS = {11B75 (11T30)},
  MRNUMBER = {2053599},
MRREVIEWER = {Ben\ Joseph\ Green},
       DOI = {10.1007/s00039-004-0451-1},
       URL = {https://doi-org.manchester.idm.oclc.org/10.1007/s00039-004-0451-1},
}

@article {Gl07,
    AUTHOR = {Glibichuk, A. A.},
     TITLE = {Properties of the sums and products of subsets in a finite
              field of prime order},
   JOURNAL = {Chebyshevski\u i\ Sb.},
  FJOURNAL = {Chebyshevski\u i\ Sbornik},
    VOLUME = {8},
      YEAR = {2007},
    NUMBER = {2(22)},
     PAGES = {30--43},
      ISSN = {2226-8383,2587-7119},
   MRCLASS = {12E20},
  MRNUMBER = {2641000},
}

@article{vD26,
AUTHOR = {{\noopsort{Doorn}}van Doorn, W. and GPT-6 Astra Pro},
TITLE = {Practical numbers and {E}gyptian fractions},
YEAR = {2026},
URL = {https://github.com/Woett/ChatGPT-s-note-on-Erdos18/blob/main/Practical%20numbers%20and%20Egyptian%20fractions.pdf},
JOURNAL = {GitHub: Woett/ChatGPT-s-note-on-Erdos18},
}

@article {Er50,
  AUTHOR = {Erd\H{o}s, P.},
  TITLE = {Az $1/x_1 + 1/x_2 + \ldots + 1/x_n  =A/B$ egyenlet eg\'{e}sz sz\'{a}m\'{u} megold\'{a}sair\'{o}l},
  JOURNAL = {Mat. Lapok},
  VOLUME = {1},
  YEAR = {1950},
  PAGES = {192-210},
}

@article {ErGr80,
    AUTHOR = {Erd\H{o}s, P. and Graham, R. L.},
     TITLE = {Old and new problems and results in combinatorial number
              theory},
    VOLUME = {28},
 PUBLISHER = {Universit\'e{} de Gen\`eve, L'Enseignement Math\'ematique,
              Geneva},
      YEAR = {1980},
     PAGES = {128},
   MRCLASS = {10-02 (05-02)},
  MRNUMBER = {592420},
MRREVIEWER = {L.\ C.\ Eggan},
}

@article{AI26,
 AUTHOR = {OpenAI},
 TITLE = {Short {E}gyptian fractions},
 YEAR = {2026},
 JOURNAL = {GitHub:openai/math},
 URL = {https://github.com/openai/math/blob/main/preprints/Short-Egyptian-fractions-September-25-2026/Short-Egyptian-fractions-September-25-2026.pdf},
}

@article{So05,
 AUTHOR = {K. Soundararajan},
 TITLE = {Approximating $1$ from below using $n$ {E}gyptian fractions},
 YEAR = {2005},
 JOURNAl = {arXiv:0502247}
}

@article {Cu22,
    AUTHOR = {Curtiss, D. R.},
     TITLE = {On {K}ellogg's {D}iophantine {P}roblem},
   JOURNAL = {Amer. Math. Monthly},
  FJOURNAL = {American Mathematical Monthly},
    VOLUME = {29},
      YEAR = {1922},
    NUMBER = {10},
     PAGES = {380--387},
      ISSN = {0002-9890,1930-0972},
   MRCLASS = {99-04},
  MRNUMBER = {1520110},
       DOI = {10.2307/2299023},
       URL = {https://doi-org.manchester.idm.oclc.org/10.2307/2299023},
}

@article {BlEl22,
    AUTHOR = {Bloom, Thomas F. and Elsholtz, Christian},
     TITLE = {Egyptian fractions},
   JOURNAL = {Nieuw Arch. Wiskd. (5)},
  FJOURNAL = {Nieuw Archief voor Wiskunde. Vijfde Serie},
    VOLUME = {23},
      YEAR = {2022},
    NUMBER = {4},
     PAGES = {237--245},
      ISSN = {0028-9825},
   MRCLASS = {11-03 (01-01 01A16 11D68)},
  MRNUMBER = {4544007},
}

@article {BlEr76,
    AUTHOR = {Bleicher, M. N. and Erd\H{o}s, P.},
     TITLE = {Denominators of {E}gyptian fractions},
   JOURNAL = {J. Number Theory},
  FJOURNAL = {Journal of Number Theory},
    VOLUME = {8},
      YEAR = {1976},
    NUMBER = {2},
     PAGES = {157--168},
      ISSN = {0022-314X,1096-1658},
   MRCLASS = {10A30},
  MRNUMBER = {404121},
MRREVIEWER = {S.\ L. G. Choi},
       DOI = {10.1016/0022-314X(76)90098-6},
       URL = {https://doi-org.manchester.idm.oclc.org/10.1016/0022-314X(76)90098-6},
}

@article {Vi47,
    AUTHOR = {Vinogradov, I. M.},
     TITLE = {The method of trigonometrical sums in the theory of numbers},
   JOURNAL = {Trav. Inst. Math. Stekloff},
  FJOURNAL = {Trav. Inst. Math. Stekloff},
    VOLUME = {23},
      YEAR = {1947},
     PAGES = {109},
   MRCLASS = {10.0X},
  MRNUMBER = {29417},
MRREVIEWER = {L.\ Schoenfeld},
}

@article{Ta21,
AUTHOR = {T. Takenouchi},
TITLE = {On an indeterminate equation},
JOURNAL = {Proceedings of the Physico-Mathematical Society of Japan, Third Series},
VOLUME = {3},
YEAR = {1921},
PAGES = {78--92},
}

@article {Ko14,
    AUTHOR = {Konyagin, S. V.},
     TITLE = {Double exponential lower bound for the number of
              representations of unity by {E}gyptian fractions},
      NOTE = {Translation of Mat. Zametki {\bf 95} (2014), no. 2, 312--316},
   JOURNAL = {Math. Notes},
  FJOURNAL = {Mathematical Notes},
    VOLUME = {95},
      YEAR = {2014},
    NUMBER = {1-2},
     PAGES = {277--281},
      ISSN = {0001-4346,1573-8876},
   MRCLASS = {11D68 (11A67 11N37)},
  MRNUMBER = {3267215},
MRREVIEWER = {Eric\ Errthum},
       DOI = {10.1134/S0001434614010295},
       URL = {https://doi-org.manchester.idm.oclc.org/10.1134/S0001434614010295},
}

@article {BaSa84,
    AUTHOR = {Balog, A. and S\'ark\"ozy, A.},
     TITLE = {On sums of integers having small prime factors. {I}, {II}},
   JOURNAL = {Studia Sci. Math. Hungar.},
  FJOURNAL = {Studia Scientiarum Mathematicarum Hungarica. Combinatorics,
              Geometry and Topology (CoGeTo)},
    VOLUME = {19},
      YEAR = {1984},
    NUMBER = {1},
     PAGES = {35--47, 81--88},
      ISSN = {0081-6906,1588-2896},
   MRCLASS = {11P55},
  MRNUMBER = {787784},
MRREVIEWER = {R.\ J.\ Cook},
}

@article{vDTa25, 
    TITLE = {The smallest denominator not contained in a unit fraction decomposition of 1 with fixed length}, 
    DOI = {10.1017/S0305004126102102}, 
    JOURNAL = {Math. Proc. Cambridge Philos. Soc.},
    FJOURNAL = {Mathematical Proceedings of the Cambridge Philosophical Society}, 
    AUTHOR = {{\noopsort{Doorn}}van Doorn, W. and Tang, Q.},
    YEAR = {2026}, 
    PAGES = {1--9}}

@article {El16,
    AUTHOR = {Elsholtz, Christian},
     TITLE = {Egyptian fractions with odd denominators},
   JOURNAL = {Q. J. Math.},
  FJOURNAL = {The Quarterly Journal of Mathematics},
    VOLUME = {67},
      YEAR = {2016},
    NUMBER = {3},
     PAGES = {425--430},
      ISSN = {0033-5606,1464-3847},
   MRCLASS = {11D68},
  MRNUMBER = {3556494},
MRREVIEWER = {Lajos\ Hajdu},
       DOI = {10.1093/qmath/haw020},
       URL = {https://doi-org.manchester.idm.oclc.org/10.1093/qmath/haw020},
}

@article{ElPl21,
  AUTHOR  = {Elsholtz, Christian and Planitzer, Stefan},
  TITLE   = {Sums of four and more unit fractions and approximate parametrizations},
  JOURNAL = {Bull. London Math. Soc.},
  FJOURNAL = {Bulletin of the London Mathematical Society},
  VOLUME  = {53},
  NUMBER  = {3},
  PAGES   = {695--709},
  YEAR    = {2021},
  DOI     = {10.1112/blms.12452}
}

@incollection {HaIo08,
    AUTHOR = {Hart, Derrick and Iosevich, Alex},
     TITLE = {Sums and products in finite fields: an integral geometric
              viewpoint},
 BOOKTITLE = {Radon transforms, geometry, and wavelets},
    SERIES = {Contemp. Math.},
    VOLUME = {464},
     PAGES = {129--135},
 PUBLISHER = {Amer. Math. Soc., Providence, RI},
      YEAR = {2008},
      ISBN = {978-0-8218-4327-7},
   MRCLASS = {11B75 (11T24 11T30)},
  MRNUMBER = {2440133},
MRREVIEWER = {Radhakrishnan\ Nair},
       DOI = {10.1090/conm/464/09080},
       URL = {https://doi-org.manchester.idm.oclc.org/10.1090/conm/464/09080},
}

@article{Sch33,
  AUTHOR   = {Schnirelmann, L.},
  TITLE    = {{\"U}ber additive Eigenschaften von Zahlen},
  FJOURNAL = {Mathematische Annalen},
  JOURNAL  = {Math. Ann.},
  VOLUME   = {107},
  PAGES    = {649--690},
  YEAR     = {1933},
  DOI      = {10.1007/BF01448914}
}

@article {Vo85,
    AUTHOR = {Vose, Michael D.},
     TITLE = {Egyptian fractions},
   JOURNAL = {Bull. London Math. Soc.},
  FJOURNAL = {The Bulletin of the London Mathematical Society},
    VOLUME = {17},
      YEAR = {1985},
    NUMBER = {1},
     PAGES = {21--24},
      ISSN = {0024-6093,1469-2120},
   MRCLASS = {11D68 (11D85)},
  MRNUMBER = {766441},
MRREVIEWER = {W.\ A.\ Webb},
       DOI = {10.1112/blms/17.1.21},
       URL = {https://doi-org.manchester.idm.oclc.org/10.1112/blms/17.1.21},
}

\end{document}